\documentclass[12pt,reqno]{amsart}
\advance \topmargin by -\headheight
\advance \topmargin by -\headsep
\evensidemargin \oddsidemargin
\newtheorem{theorem}{Theorem}[section]
\newtheorem{lemma}[theorem]{Lemma}

\usepackage{amsmath, amsthm, amssymb, amsfonts}
 \usepackage{hyperref}
\usepackage{enumerate}
\usepackage{url}
\usepackage{dsfont}
\usepackage{mathrsfs}

\usepackage{bm}

\usepackage{xcolor}

\usepackage{fancyvrb}

\theoremstyle{definition}

\theoremstyle{remark}

\numberwithin{equation}{section}

\DeclareMathOperator{\Mod}{mod}
\newcommand{\mmod}[1]{\;(\Mod{ #1})}

\renewcommand{\geq}{\geqslant}
\renewcommand{\leq}{\leqslant}

\def\eps{\varepsilon}

\def \cB {\mathcal B}

\def \cD {\mathcal D}

\def \Z {\mathbb{Z}}

\newcommand{\n}[1]{\|{#1}\|}
\newcommand\1{\mathds{1}}

\title{Counting $2\times 2$ integer matrices with a given determinant}

\author{Jonathan Chapman \and Akshat Mudgal}
\address{Mathematics Institute, Zeeman Building, University of Warwick, Coventry CV4 7AL, United Kingdom}
\email{Jonathan.Chapman@warwick.ac.uk}
\address{Mathematics Institute, Zeeman Building, University of Warwick, Coventry CV4 7AL, United Kingdom}
\email{Akshat.Mudgal@warwick.ac.uk}

 \subjclass[2020]{11D45 (primary); 11D09, 11N37 (secondary)} 
\keywords{Integer matrices with fixed determinant, Restricted divisor correlations}

\begin{document}

\begin{abstract}
Given positive integers $h, N$ satisfying $1 \leqslant h \leqslant 2N^2$, we define $T(h,N)$ to be the number of $2\times 2$ integer matrices with determinant equal to $h$ whose entries lie in $[-N,N]$. 
Our main result states that for any $\varepsilon >0$, one has 
\[ T(h,N) = \frac{16}{\zeta(2)} N^2 \bigg( \sum_{d |h} \frac{1}{d}  \bigg) + O_{\varepsilon}(N^{\varepsilon} (N+ h)).\]
This quantitatively improves upon recent work of Afifurrahman and Ganguly--Guria, and delivers square-root cancellation estimates when $h \leq N$. We further show that when $h$ is large, the error term is of approximately the correct order. 
\end{abstract}
\maketitle

\section{Introduction}

For any integer $h$, let 
\[ \mathcal{D}_h = \{ (a,b,c,d) \in \mathbb{Z}^4 : ad - bc = h\}.\]
One can naturally interpret $\mathcal{D}_h$ as the set of $2\times 2$ integer matrices with determinant $h$.
There are various classical questions in analytic number theory that concern counting elements of $\cD_h$ in expanding regions. For instance, writing $\cB(N)$ to be the ball of radius $N$ in $\mathbb{R}^4$ for every $N \in \mathbb{N}$, we can define 
\[ T_{\ell^2}(h,N) = |\cD_h \cap \cB(N)|  =  6 N^2 (\sum_{d|h} 1/d)  + E_{\ell^2}(h,N) . \]
 A classical result of Selberg states that 
 \begin{equation} \label{selbergresult}
 E_{\ell^2}(1,N) \ll N^{4/3}.
 \end{equation}
 This bound has not been improved, and in fact it is conjectured that
 \begin{equation} \label{conjl2}
     E_{\ell^2}(1,N) \ll_{\eps} N^{1  + \varepsilon},
 \end{equation}
 see \cite[Chapter 12]{Iw2002}. A more modern version of \eqref{selbergresult} is recorded by Iwaniec \cite[Theorem 12.4]{Iw2002}, which implies that 
$E_{\ell^2}(h,N)\ll h^{1/3} N^{4/3} \sum_{d|h}1/d$ for all $1 \leq h \leq N^2$.

One can also consider counting $(a,b,c,d) \in \cD_h$ with $a,b,c,d \in \mathbb{N}$ and $bc \leq N^2$. Thus, we define
\[  T_{\rm div}(h, N) = \sum_{1 \leq n \leq N^2} d(n) d(n+h) = M_{\rm div}(h,N) + E_{\rm div}(h,N), \]
where $d(n) = \sum_{x,y \in \mathbb{N}} \1_{n = xy}$ counts the number of divisors of $n$ and $M_{\rm div}(h,N)$ denotes the main term described in \cite[(1.2)]{Mo1994}. The problem of finding good quantitative estimates for the above is known as the binary additive divisor problem, see \cite{DI1982, Es1930, IM1995, Meu2001, Mo1994} and the references therein. This is a very well-studied problem, in part, due to its close connections to the fourth moment of the Riemann zeta function, see \cite{HB1979}. The work of Motohashi \cite{Mo1994}, which relies on spectral methods, implies that for any $h \leq N^{40/27}$  one has
\[ E_{\rm div}(h,N) \ll_{\eps} N^{4/3 + \eps},\]
see also \cite[Theorem 1]{Mo1994} for a modified asymptotic formula which holds for a wider range of $h$. Motohashi \cite[Theorem 7]{Mo1994} also proved that for fixed $h$ and large $N$, one has
\begin{equation} \label{dnd4}
E_{\rm div}(h,N) \gg N.
\end{equation}
It is natural to speculate that this should be almost optimal. In fact, Conrey--Gonek \cite[Conjecture 3]{CG2001} have conjectured that for all $0 < h \ll N$, one should have
\begin{equation} \label{divconj}
    E_{\rm div}(h,N) \ll_{\eps} N^{1 + \eps}.
\end{equation}


From the perspective of counting integer solutions to quadratic equations, a natural question concerns studying points in $\cD_h \cap [-N,N]^4$. This type of problem arises in various contexts including combinatorial geometry \cite{MWY2022}, counting commuting matrices \cite{CM2025, Mu2024}, diophantine approximation \cite{DS1941, BCGW2018}, theory of random multiplicative functions \cite{HNR2015, HL2016}, and  counting integer solutions to quadratic equations \cite{HB1996, Mc2024}. Thus, we define
\[ T(h,N)  = |\cD_h \cap [-N,N]^4| = \sum_{n \in \mathbb{Z}} d'(n) d'(n+h),  \ \ \text{where} \ \ d'(n) = \sum_{|x|,|y| \leq N} \1_{n= xy} \]
is the restricted divisor function. We further write
\begin{equation} \label{defthn}
 T(h,N) = \frac{16}{\zeta(2)} N^2 \big( \sum_{d|h}1/d \big) + E_{\rm \ell^{\infty}}(h,N) . 
\end{equation}
  This can be construed as a circle-method heuristic; indeed, the first term here can be written as a product of a singular series and a singular integral when $h \ll_{\eps} N^{2 - \eps}$, see \cite{CM2026}.

The problem of estimating $T(h,N)$ was first analysed by Afifurrahman \cite{Af2024}, who used results on counting points on modular hyperbola, which themselves require deep estimates of Weil for Kloosterman sums, to prove that for every $1 \leq h \leq 2N^2$ one has
\begin{equation} \label{res1}
E_{\rm \ell^{\infty}}(h,N)  \ll_{\varepsilon}N^{\varepsilon}( h + N^{5/3}).
\end{equation}
Subsequently, Ganguly--Guria \cite{GG2024} used Fourier analysis and spectral methods to show that in the shorter range  $1 \leq h \leq N^{1/3}$ one has 
\begin{equation} \label{res2}
    E_{\rm \ell^{\infty}}(h,N)  \ll_{\varepsilon} h^{\theta} N^{3/2 + \varepsilon},
\end{equation}
where $\theta$ is any admissible exponent for the generalised Ramanujan conjecture. Deep work of Kim--Sarnak \cite{Kim2003} implies that $\theta \leq 7/64$.

Noting the aforementioned bounds for $E_{\ell^2}(h,N)$ and $E_{\rm div}(h,N)$, it is reasonable to speculate that the bounds in \eqref{res1} and \eqref{res2} can be improved significantly. We are able to do this in a short and completely elementary fashion.

\begin{theorem} \label{thm1.1}
    For any $h,N \in \mathbb{N}$ with $h \leq 2N^2$, one has 
\[ T(h,N) = \frac{16}{\zeta(2)} N^2 \big( \sum_{d |h} 1/d  \big) + O_{\varepsilon}(N^{\varepsilon} (h + N)). \]
\end{theorem} 

This improves upon \eqref{res1} and \eqref{res2}, while also confirming an analogue of the  conjectured bounds \eqref{conjl2} and \eqref{divconj} in our setting. Indeed, Theorem \ref{thm1.1} implies that for all $0 < h \leq N$ one has
\begin{equation} \label{errortermepsilon}
E_{\ell^{\infty}}(h,N) \ll_{\eps} N^{1 + \eps}.
\end{equation}  
We note that our proof can be analysed more carefully to give a more precise characterisation of the $N^{\eps}$ factor in the error term, but we have not pursued this here. 

Our methods also deliver an analogue of Motohashi's lower bound \eqref{dnd4} for $E_{\ell^{\infty}}(h,N)$ for much larger shifts $h$. 

\begin{theorem} \label{prop1.5}
    Let $\delta \in (0,1]$, let $N \in \mathbb{N}$ be sufficiently large in terms of $\delta$ and let $h \in \mathbb{N}$ satisfy $N^{1 + \delta} \leq h \leq 2N^2$. Then
    \[ E_{\ell^{\infty}}(h,N) \gg h.\]
\end{theorem}

Thus, whenever $h \geq N^{1+ \delta}$ for any fixed $\delta>0$, one cannot expect \eqref{errortermepsilon} to hold, and the error term in Theorem \ref{thm1.1} turns out to be roughly of the right order.

An interesting aspect of this lower bound is that it holds in the range $h \gg N^2$, which in turn means that the error term in \eqref{defthn} matches the order of the main term for many choices of $h \gg N^2$. This does not mean that the circle-method heuristic for \eqref{defthn} fails completely when $h \gg N^2$, it just means that one expects a different main term because the size of the singular integral varies as $h$ varies, see \cite{CM2026}. In fact, the latter setting is much harder than the regime $h = o(N^2)$ studied in this paper and nothing seems to have been previously known for $T(h,N)$ when $N^2 \ll  h \leq (2 - o(1))N^2$. Moreover, the only known result for $T_{\ell^2}(h,N)$ in this setting seems to follow from a general class of results proven by Oh \cite{Oh2004} using techniques from dynamics.  While Theorem \ref{prop1.5} implies that our methods in this paper cannot  yield non-trivial estimates for $T(h,N)$ when $h \gg N^2$, we use completely different methods to analyse this setting in \cite{CM2026}.

We remark that obtaining an asymptotic formula in the case when $h=0$ is significantly simpler than the cases when $h \neq 0$ since the former allows
for many additional symmetries. Indeed, writing $T(0,N)$ to be the number of solutions to $ad =bc$ with $a,b,c,d \in [-N,N]$, one can swiftly show that 
\begin{equation} \label{singularmatrixcount} T(0,N) = \frac{16}{\zeta(2)}N^2 \log N + O(N^2). \end{equation}
An asymptotic formula with an explicit second order term of order $N^2$ can be deduced from work of Ayyad--Cochrane--Zheng, see \cite[Theorem 3]{ACZ1996}. 


As previously mentioned, our results have a natural interpretation as counting integer matrices with fixed determinant and bounded $\ell^{\infty}$ norm. Similarly, $T_{\ell^2}(1,N)$ counts matrices $\gamma \in {\rm SL}_2(\mathbb{Z})$ such that $\n{\gamma} = {\rm Tr}(\gamma^T \gamma)^{1/2} \leq N$. The latter type of results have been generalised to matrices $\gamma \in {\rm SL}_n(\mathbb{Z})$ with $n \geq 3$, see, for instance, the very nice work of Duke--Rudnick--Sarnak \cite{DRS1993}, and subsequent results by Gorodnik--Nevo--Yehoshua \cite{GNY2017} and Blomer--Lutsko \cite{BL2024}. A higher dimensional analogue of \eqref{singularmatrixcount} in the bounded $\ell^2$ case was established by Katznelson \cite{Ka1993} via techniques from the geometry of numbers. It would also be interesting to consider versions of these higher dimensional results in the bounded $\ell^{\infty}$ case.

\subsection*{Outline} We use \S\ref{sec2} to record various standard lemmas which we will employ in our proof of Theorems \ref{thm1.1}.  We begin our proof of Theorem \ref{thm1.1} in \S\ref{sec3} by performing some preliminary manoeuvres to reduce our problem to analysing a weighted count of points in arithmetic progressions with varying lengths and moduli. We are able to make these arithmetic progressions slightly more uniform in size at the cost of introducing the error term appearing in Theorem \ref{thm1.1}. In \S\ref{sec4}, we use this uniformisation trick to further reduce our problem to estimating pairs of coprime integers in large intervals and counting coprime residue classes satisfying extra congruence conditions, both of these being standard results from elementary number theory.  We conclude \S\ref{sec4} by proving Theorem \ref{prop1.5}.  

\subsection*{Notation} We employ Vinogradov notation, that is, we write $Y \ll_{z} X$, or equivalently $Y =O_z(X)$, to mean that $|Y| \leq C_z X$, where $C_z>0$ is some constant depending on the parameter $z$. Unless stated otherwise, whenever $\eps$ appears in any bound, it will mean that the bound holds for every $\eps >0$, though the implicit constant may depend on $\eps$. We denote the greatest common divisor of two integers $a$ and $b$ by $(a,b)$.  

\subsection*{Acknowledgements} We thank Sam Chow, Lasse Grimmelt, Jori Merikoski, V. Vinay Kumaraswamy, and Trevor Wooley for helpful comments. JC is supported by EPSRC through Joel Moreira's Frontier Research Guarantee grant, ref. \texttt{EP/Y014030/1}. AM is supported by a Leverhulme Early Career Fellowship \texttt{ECF-2025-148}.

Finally, as we were finishing the first version of this paper, it came to our attention that Dhanda--Haynes--Prasala \cite{DHP2025} have independently proved Theorem \ref{thm1.1}.  

\subsection*{Rights}

For the purpose of open access, the authors have applied a Creative Commons Attribution (CC-BY) licence to any Author Accepted Manuscript version arising from this submission.

\section{Auxiliary lemmas} \label{sec2}

We utilise this section to record various standard results from elementary number theory that we will use throughout our paper. 

\begin{lemma} \label{lem2.1}
    Let $1 \leq q \leq u$ and $r$ be positive integers such that $q|u$ and $(r, q) = 1$. Then 
    \[ \sum_{\substack{1\leq v < u, \\ (v,u) = 1}} \1_{v \equiv r \mmod{q}} = \frac{\varphi(u)}{ \varphi(q)}. \]
\end{lemma}

\begin{proof}
 Let $F:(\Z/u\Z)^\times\to(\Z/q\Z)^\times$ be the homomorphism $F:t+u\Z\mapsto t+q\Z$. Since $F$ is surjective, the first isomorphism theorem shows that
    \begin{equation*}
        |\{v\in(\Z/u\Z)^\times : v \equiv r \mmod{q}\}| = |F^{-1}(\{r\})| = \frac{|(\Z/u\Z)^\times|}{|(\Z/q\Z)^\times|} = \frac{\varphi(u)}{\varphi(q)}. \qedhere
    \end{equation*}
\end{proof}

\begin{lemma} \label{lem2.2}
  Let $m,y,z$ be positive integers. Then
  \[ \sum_{t \in (\mathbb{Z}/y \mathbb{Z})^{\times} } \1_{tm \equiv z \mmod{y}} \leq (y,m). \]
\end{lemma}

\begin{proof}
Any $t \in (\mathbb{Z}/y \mathbb{Z})^{\times}$ satisfying $t m \equiv z \mmod{y}$ further satisfies
\[ t (m/(y,m)) \equiv z/(y,m) \mmod{y/(y,m)}. \]
Note that the left-hand side of the above congruence is coprime to $y/(y,m)$, and so, for there to be any solutions to this congruence, $z/(y,m)$ must be coprime to $y/(y,m)$. In this case, we can apply Lemma \ref{lem2.1} to obtain
\[ \sum_{t \in (\mathbb{Z}/y \mathbb{Z})^{\times} } \1_{tm \equiv z \mmod{y}} \leq \sum_{\substack{1 \leq t < y, \\ (t,y) = 1}} \1_{t \equiv \frac{z}{(y,m)} (\frac{m}{(y,m)})^{-1} \mmod{\frac{y}{(y,m)}}} = \frac{\varphi(y)}{\varphi(y/(y,m))} \leq (y,m). \qedhere \]
\end{proof}

\begin{lemma} \label{lem2.3}
    Let $m$ be a positive integer and let $M \geq 1$ be a real number. Then 
    \[ \sum_{1 \leq y \leq M} (y,m) \ll_{\varepsilon} M m^{\varepsilon} \]
\end{lemma}

\begin{proof}
    Using the divisor bound, we have
    \[ \sum_{1 \leq y \leq M} (y,m) \leq \sum_{\substack{1 \leq k \leq M, \\ k |m}} k \sum_{\substack{k \leq y \leq M,\\ k |y }} 1 \leq \sum_{\substack{1 \leq k \leq M, \\ k |m}} k (M/k) \ll_{\varepsilon} M m^{\varepsilon}. \qedhere \]
\end{proof}

We will also require the following standard estimate on the number of elements of an interval which are coprime to some fixed moduli, see  \cite[Theorem 8.29]{NZM1991}.

\begin{lemma} \label{lem2.4}
    Let $q,Y$ be positive integers and let $X \geq 1$ be a real number. Then 
    \[ \sum_{Y\leq n < Y + X} \mathds{1}_{(n,q) = 1}  = \frac{\varphi(q)}{q} X + O(\tau'(q) ),\]
    where $\tau'(q)=\sum_{d\mid q}|\mu(d)|$ counts the number of square-free divisors of $q$.
\end{lemma}

\section{Preliminary manoeuvres} \label{sec3}

We employ this section and the next to present our proof of Theorem \ref{thm1.1}. Throughout this section, let $N$ be some natural number and let $1 \leq h \leq 2N^2$. Upon excluding solutions to $ax - by = h$ with $abxy = 0$ and noting the various symmetries, we see that
\[ T(h,N) = 4 \sum_{\substack{d|h, \\ d \leq N}}  \big(  \sum_{\substack{1 \leq u,v \leq N/d \\ (u,v) = 1}} \ \sum_{x,y \in [-N,N]} \mathds{1}_{ux - vy = h/d} \big)  + O_{\varepsilon}(N^{1 + \varepsilon}). \]
 Excluding further the contribution of the solutions when $u = v = 1$, which contribute at most $O_{\varepsilon}(N^{1 + \varepsilon})$ due to the divisor bound, and then noting that in all other cases we may assume without loss of generality that $u >v$ (as $x,y \in [-N,N]$), we obtain
\begin{equation} \label{eqn3.1}
    T(h,N)  = 8 \sum_{\substack{d|h, \\ d \leq N}}  \big(  \sum_{\substack{1 \leq v < u \leq N/d \\ (u,v) = 1}} \ \sum_{x,y \in [-N,N]} \mathds{1}_{ux - vy = h/d} \big)  + O_{\varepsilon}(N^{1+\varepsilon}). 
\end{equation}
Hence, given $n\in \mathbb{N}$ and coprime integers $1\leq v < u$, we define
\[ r_{u,v}(n) = \sum_{x,y \in [-N,N]} \1_{ux - vy = n} \ \ \text{and} \ \  \tilde{r}_{u,v}(n) = \sum_{y\in [-N,N]} \1_{y \equiv -n v^{-1} \mmod{u} }. \]
Note that $y \in [-N,N]$ can lead to a valid solution of $ux - vy = n$ with $x \in [-N,N]$ if and only if 
\[ y \equiv -n v^{-1} \mmod{u} \ \ \text{and} \ \ y \in [-(Nu+n)/v, (Nu-n)/v] \cap [-N,N].\]
Furthermore, since $u/v >1$ and $n \geq 1$, we always have $-(Nu+n)/v < - N$. This immediately implies that 
\begin{align} \label{eqn3.2}
    0 \leq \tilde{r}_{u,v}(n) - r_{u,v}(n)  & \leq \sum_{ (Nu -n)/v <  y \leq N} \1_{y \equiv - n v^{-1} \mmod{u}} \nonumber \\
    & = \bigg(\frac{n - N(u-v)}{uv} + O(1) \bigg) \1_{n \geq N(u-v)} . 
\end{align}
In fact, if $1 \leq n < N(u-v)$, then $r_{u,v}(n) = \tilde{r}_{u,v}(n)$, and if $n > N(u+v)$, then $r_{u,v}(n) = 0$ while $|\tilde{r}_{u,v}(n) -2N/u| \leq 5$. Moreover, when $N(u-v) \leq n \leq N(u+v)$, equality holds in the second inequality in \eqref{eqn3.2}.

Setting
\[ \tilde{T}(h,N) =   8 \sum_{\substack{d|h, \\ d \leq N}}  \big(  \sum_{\substack{1 \leq v < u \leq N/d \\ (u,v) = 1}} \ \tilde{r}_{u,v}(h/d) \big) ,\]
our proof of Theorem \ref{thm1.1} now reduces to proving that $\tilde{T}(h,N)$ is a good approximation for $T(h,N)$ and then appropriately estimating $\tilde{T}(h,N)$. The following lemma proves the former statement.

\begin{lemma} \label{lem3.1}
    Let $h,N$ be natural numbers with $1 \leq h \leq 2N^2$, let $\varepsilon>0$. Then 
    \[ |T(h,N) - \tilde{T}(h,N)| \ll_{\varepsilon} N^{\varepsilon}(h + N).  \]
\end{lemma}

\begin{proof}
Noting \eqref{eqn3.1} and \eqref{eqn3.2}, we have
\begin{align}  \label{eqn3.3}
\tilde{T}(h,N) - & T(h,N)
 \ll \sum_{\substack{d|h, \\ d \leq N}}\sum_{\substack{ 1 \leq v < u \leq N/d, \\ (u,v) =1}} \bigg(\frac{h/d- N(u-v)}{uv} + O(1) \bigg) \1_{h/d \geq N(u-v)}  + O_{\varepsilon}(N^{1 + \varepsilon}) \nonumber \\
& =  \sum_{\substack{d|h, \\ d \leq N}} \sum_{\substack{ \max\{1,u - h/(dN)\} \leq v < u \leq N/d, \\ (u,v) =1}}  \bigg(\frac{h}{duv} - \frac{N(u-v)}{uv} + O(1) \bigg) + O_{\varepsilon}(N^{1 + \varepsilon}) .
\end{align}
Note that
\[\sum_{\substack{d|h, \\ d \leq N}}\sum_{\substack{ \max\{1,u - h/(dN)\}\leq v < u \leq N/d, \\ (u,v) =1}} \frac{h}{duv} \ll  \sum_{d|h}  \frac{h}{d} (\log N)^2 \ll_{\varepsilon} h N^{\varepsilon} \]
and
\[     \sum_{\substack{d|h, \\ d \leq N}} \sum_{\substack{ \max\{1,u - h/(dN)\} \leq v < u \leq N/d, \\ (u,v) =1}} \frac{N(u-v)}{uv} \leq  \sum_{d|h} \frac{h}{d} \sum_{1 \leq v \leq u \leq N/d}\frac{1}{uv} \ll_{\varepsilon} h N^{\varepsilon},\]
where in the second inequality, we have used the fact that $u-v \leq h/(dN)$. Finally, we see that
\begin{align*}
    \sum_{\substack{d|h, \\ d \leq N}} \ \sum_{\substack{ \max\{1,u - h/(dN)\} \leq v < u \leq N/d, \\ (u,v) =1}}1 
    & \ll \sum_{d|h} \sum_{1 \leq u \leq N/d} \bigg(\frac{h}{dN} +1\bigg) \ll \sum_{d|h}  \frac{h}{dN}\bigg( \frac{N}{d} + O(1)\bigg) + O_{\varepsilon}(N^{1 + \varepsilon}) \\
    & \ll_{\varepsilon} N^{\varepsilon}(h + N).
\end{align*}
The preceding three inequalities combine with \eqref{eqn3.3} to deliver the claimed estimate. 
\end{proof}

Thus, it suffices to prove Theorem \ref{thm1.1} for $\tilde{T}(h,N)$, and we will do so in the next section.

\section{Proof of Theorems \ref{thm1.1} and  \ref{prop1.5}} \label{sec4}

Our first aim of this section is to prove the following result. 

\begin{lemma} \label{lem4.1}
    Let $h,N$ be natural numbers with $1 \leq h \leq 2N^2$. Then 
    \[ \tilde{T}(h,N) = \frac{16}{\zeta(2)} N^2 \big(\sum_{d |h} 1/d \big) + O_{\varepsilon}(N^{1+\varepsilon}).  \]
\end{lemma}

\begin{proof}
Let $1 \leq d \leq h$ and $1 \leq v < u$ be integers such that $(u,v) = 1$ and $d|h$, and suppose that $y \equiv -(h/d) v^{-1} \mmod{u}$. Writing $k = (h/d,u)$, we see that this is equivalent to the congruence $(y/k) \equiv -(h/(dk)) v^{-1} \mmod{u/k}$. Moreover, writing $y' = y/k$, we see that $y'$ is coprime to $u/k$ since  $(v,u/k) = (h/(dk), u/k) = 1$. Thus, we have that
\[ \sum_{\substack{1 \leq v < u \leq N/d \\ (u,v) = 1}} \ \tilde{r}_{u,v}(h/d) = \sum_{1 < u \leq N/d} \  \sum_{\substack{y' \in [-N/k, N/k], \\ (y', u/k) = 1 }} \  \sum_{\substack{1 \leq v < u, \\ (v,u) = 1}} \1_{y'  \equiv -(h/(dk)) v^{-1} \mmod{u/k}}. \]
Applying Lemma \ref{lem2.1} and Lemma \ref{lem2.4}, we find that
\begin{align*}
\sum_{\substack{1 \leq v < u \leq N/d \\ (u,v) = 1}} \ \tilde{r}_{u,v}(h/d) 
& = \sum_{1 < u \leq N/d} \  \sum_{\substack{y' \in [-N/k,N/k], \\ (y', u/k) = 1 }} \  \frac{\varphi(u)}{\varphi(u/k)}   \\
& = \sum_{1 \leq u \leq N/d} \frac{\varphi(u)}{\varphi(u/k)}  \bigg( \frac{2N}{k} \frac{\varphi(u/k)}{u/k} + O(\tau'(u/k)) \bigg) \\
& = 2N \sum_{1 \leq u \leq N/d} \frac{\varphi(u)}{u} + O \bigg( \sum_{1 \leq u \leq N/d} \frac{\varphi(u)}{\varphi(u/k)} \tau'(u/k) \bigg).
\end{align*}
Hence,
\begin{align} \label{eqn4.1}
\tilde{T}(h,N) & = 16 N \sum_{\substack{d|h, \\ d \leq N}} \sum_{1 \leq u \leq N/d}  \frac{\varphi(u)}{u} + O \bigg(\sum_{\substack{d|h, \\ d \leq N}}\sum_{1 \leq u \leq N/d} \frac{\varphi(u)}{\varphi(u/k)} \tau'(u/k) \bigg) \nonumber \\
& = \frac{16 N^2}{\zeta(2)} \big( \sum_{\substack{d|h, \\ d \leq N}}1/d \big) + O \bigg(\sum_{d|h} \sum_{1 \leq u \leq N/d} \frac{\varphi(u)}{\varphi(u/k)} \tau'(u/k) \bigg) +  O_{\varepsilon}(N^{1 + \varepsilon}).
\end{align}
Since $k = (u, h/d)$ and
\[ \frac{\varphi(u)}{\varphi(u/k)} = k \prod_{\substack{\text{primes }p, \ p | u \\ p \nmid (u/k)}} (1 - 1/p) \leq k , \]
we find that the error term in \eqref{eqn4.1} is
\begin{align*}
& \ll \sum_{d|h} \sum_{1 \leq u \leq X/d} (u, h/d) \tau'(u/k)  \ll_{\varepsilon} N^{\varepsilon/3} \sum_{d|h} \sum_{1 \leq u \leq N/d} (u, h/d) \ll_{\varepsilon} N^{2\varepsilon/3} \sum_{d |h} N/d \ll_{\varepsilon} N^{1 + \varepsilon}
\end{align*}
with the penultimate step following from Lemma \ref{lem2.3}.
Moreover, the main term in \eqref{eqn4.1} equals
\[ \frac{16 N^2}{\zeta(2)} \big( \sum_{\substack{d|h, \\ d \leq N}}1/d \big) = \frac{16 N^2}{\zeta(2)} \big( \sum_{d|h}1/d \big) + O_{\varepsilon}(N^{1 + \varepsilon}), \]
thus concluding the proof.
\end{proof}

We remark that Lemma \ref{lem4.1} combines with Lemma \ref{lem3.1} to dispense the conclusion of Theorem \ref{thm1.1}. We now present the proof of Theorem \ref{prop1.5}.

\begin{proof}[Proof of Theorem \ref{prop1.5}]
Recalling \eqref{eqn3.1} and the definition of $\tilde{T}(h,N)$, we see that
\[ \tilde{T}(h,N) - T(h,N) =  \sum_{\substack{d|h, \\ d \leq N}}\sum_{\substack{ 1 \leq v < u \leq N/d, \\ (u,v) =1}} ( \tilde{r}_{u,v}(h/d)  - r_{u,v}(h/d) ) + O_{\varepsilon}(N^{1 + \varepsilon}).\]
Moreover, we have $\tilde{r}_{u,v}(h/d)  \geq  r_{u,v}(h/d)$ for all valid choices of $u,v,d,h$. Noting Lemma \ref{lem4.1}, our main aim will be to show that the right-hand side above is $\gg h$. In order to obtain this lower bound, we just consider the contribution from some special choices of $d,u,v$. Noting the remark following \eqref{eqn3.2}, we see that when $d=1$ and $1 \leq v < u \leq N$ satisfy $(u,v) = 1$ and $u < N/100$ and $N(u+v) < h$, then $r_{u,v}(h) = 0$ while $\tilde{r}_{u,v}(h) \geq 2N/u - 5$. 
Thus, if $d=1$ and $u < h/(100N)$, then we obtain all the above size constraints. Hence, in this case, we have $\tilde{r}_{u,v}(h) \geq 2N/u - 5  \geq  N/u$. Therefore, we deduce that
\begin{align*}
    \tilde{T}(h,N) - T(h,N) & \geq \sum_{1 < u < h/(100N) }  \sum_{1 \leq v \leq u, (v,u) = 1}  N/u + O_{\varepsilon}(N^{1 + \varepsilon}) \\
    & \geq N  \sum_{1 < u < h/(100N) } \varphi(u)/u + O_{\varepsilon}(N^{1 + \varepsilon})  .  
\end{align*}
Note that for any sufficiently large $L \in \mathbb{N}$, we have
\begin{align*}
\sum_{1 <u < L} \varphi(u)/u 
& = \sum_{1< u < L} \sum_{d|u} \mu(d)/d = \sum_{1 \leq d < L}\mu(d)/d ( L/d + O(1)) \\
& = L \sum_{1 \leq d < L} \mu(d)/d^2 + O( \log L)  = L/\zeta(2) + O(\log L) \gg L.
\end{align*}
Since $h/100N = N^{\delta}/100$, we therefore see that whenever $N$ is sufficiently large in terms of $\delta$ and $\varepsilon$, we have
\[ \tilde{T}(h,N) - T(h,N) \gg h + O_{\varepsilon}(N^{1 + \varepsilon}).  \]
Noting Lemma \ref{lem4.1} and choosing $\varepsilon$ to be sufficiently small in terms of $\delta$, we conclude that
\[ |T(h,N) - \frac{16}{\zeta(2)} N^2 \big(\sum_{d |h} 1/d \big) | \gg h. \qedhere\]
\end{proof}


\begin{thebibliography}{50}

\bibitem{Af2024}
M. Afifurrahman, \emph{A uniform formula on the number of integer matrices with given determinant and height}, J. Number Theory \textbf{281} (2026), 741--770.



\bibitem{ACZ1996}
A. Ayyad, T. Cochrane, Z. Zheng, \emph{The congruence $x_1x_2 \equiv x_3x_4 \mmod{p}$, the equation $x_1 x_2=x_3 x_4$, and mean values of character sums}, J. Number Theory \textbf{59} (1996), no. 2, 398--413.

\bibitem{BL2024}
V. Blomer, C. Lutsko, \emph{Hyperbolic lattice point counting in unbounded rank}, J. Reine Angew. Math. \textbf{812} (2024), 257–274.


\bibitem{CM2025}
J. Chapman, A. Mudgal, \emph{On commuting integer matrices}, arXiv:2504.15839, to appear in Trans. Amer. Math. Soc. 


\bibitem{CM2026}
J. Chapman, A. Mudgal, \emph{Counting solutions to the quadratic determinant equation}, arXiv:2605.15434.


\bibitem{BCGW2018}
T. F. Bloom, S. Chow, A. Gafni, A. Walker, \emph{Additive energy and the metric Poissonian property}, Mathematika 64 (2018), no. 3, 679--700.


\bibitem{CG2001}
J. B. Conrey, S. M. Gonek, \emph{High moments of the Riemann zeta-function}, Duke Math. J. \textbf{107} (2001), no. 3, 577--604.



\bibitem{DI1982}
J.-M. Deshouillers, H. Iwaniec, \emph{An additive divisor problem}, J. London Math. Soc. (2) \textbf{26} (1982), no. 1, 1--14.



\bibitem{DHP2025}
K. Dhanda, A. Haynes, S. Prasala, \emph{Counting $2\times 2$ matrices with fixed determinant and bounded coefficients}, arXiv:2509.16890.


\bibitem{DS1941}
R. J. Duffin,  A. C. Schaeffer, \emph{Khintchine's problem in metric Diophantine approximation}, Duke Math. J. \textbf{8} (1941), 243--255.



\bibitem{DRS1993}
W. Duke, Z. Rudnick, P. Sarnak, \emph{Density of integer points on affine homogeneous varieties}, Duke Math. J. \textbf{71} (1993) 143--179.


\bibitem{Es1930}
T. Estermann, \emph{On the Representations of a Number as the Sum of Two Products}, Proc. London Math. Soc. (2) \textbf{31} (1930), no. 2, 123--133.


\bibitem{GG2024}
S. Ganguly, R. Guria, \emph{Lattice points on determinant surfaces and the spectrum of the automorphic Laplacian}, arXiv:2410.04637.


\bibitem{GNY2017}
A. Gorodnik, A. Nevo, G. Yehoshua, \emph{Counting lattice points in norm balls on higher rank simple Lie groups}, Math. Res. Lett. 24 (2017), no. 5, 1285--1306.


\bibitem{HNR2015}
A. J. Harper, A. Nikeghbali, M. Radziwiłł, \emph{A note on Helson's conjecture on moments of random multiplicative functions}, Analytic number theory, 145–169, Springer, Cham, 2015.


\bibitem{HL2016}
W. P. Heap, S. Lindqvist, \emph{Moments of random multiplicative functions and truncated characteristic polynomials}, Q. J. Math. \textbf{67} (2016), no. 4, 683--714.



\bibitem{HB1979}
D. R. Heath-Brown, \emph{The fourth power moment of the Riemann zeta function}, Proc. London Math. Soc. (3)   \textbf{38} (1979), no. 3, 385--422.




\bibitem{HB1996}
D. R. Heath-Brown, \emph{A new form of the circle method, and its application to quadratic forms}, J. Reine Angew.
Math. \textbf{481} (1996), 149--206.



\bibitem{IM1995}
A. Ivi\'{c}, Y. Motohashi, \emph{On some estimates involving the binary additive divisor problem}, Quart. J. Math. Oxford Ser. (2) \textbf{46} (1995), no. 184, 471--483.



\bibitem{Iw2002}
H. Iwaniec, \emph{Spectral methods of automorphic forms}, Second edition Grad. Stud. Math., 53 American Mathematical Society, Providence, RI; Revista Matemática Iberoamericana, Madrid, 2002. xii+220 pp.
ISBN:0-8218-3160-7

\bibitem{Ka1993}
Y. R. Katznelson, \emph{Singular matrices and a uniform bound for congruence groups of ${\rm SL}_n(\mathbb{Z})$}, Duke Math. J. \textbf{69} (1993), no. 1, 121--136.


\bibitem{Kim2003}
H. Kim, \emph{Functoriality for the exterior square of ${\rm GL}_4$ and the symmetric fourth of ${\rm GL}_2$},
With appendix 1 by D. Ramakrishnan and appendix 2 by Kim and P. Sarnak, J. Amer. Math. Soc. \textbf{16} (2003), no. 1, 139--183.


\bibitem{MWY2022}
G. Martin, E. P. White, C. H. Yip, \emph{Asymptotics for the number of directions determined by $[n] \times [n]$ in $\mathbb{F}_p^2$},
Mathematika \textbf{68} (2022), no. 2, 511--534.


\bibitem{Mc2024}
O. McGrath, \emph{On the asymmetric additive energy of polynomials}, Trans. Amer. Math. Soc. 377 (2024), no. 7, 4895--4930.

\bibitem{Meu2001}
 T. Meurman, \emph{On the binary additive divisor problem}, Number theory (Turku, 1999), 223–-246,
Walter de Gruyter \& Co., Berlin, 2001.


\bibitem{Mo1994}
Y. Motohashi, \emph{The binary additive divisor problem}, Ann. Sci. \'{E}cole Norm. Sup. (4)   \textbf{27} (1994), no. 5, 529--572.



\bibitem{Mu2024}
A. Mudgal, \emph{On commuting pairs in arbitrary sets of $2 \times 2$ matrices}, arXiv:2411.10404.

\bibitem{NZM1991}
I. Niven, H.S. Zuckerman, and H.L. Montgomery, \emph{An introduction to the theory of numbers}, Fifth Edition, John Wiley \& Sons, Inc., New York, 1991.

\bibitem{Oh2004}
H. Oh, \emph{Hardy-Littlewood system and representations of integers by an invariant polynomial}, Geom. Funct. Anal. \textbf{14} (2004), no. 4, 791--809.

\end{thebibliography}
\end{document}